\documentclass[aip, sd, amsmath, amssymb, reprint]{revtex4-1}

\usepackage{amsmath}
\usepackage{amsfonts}
\usepackage{amssymb}
\usepackage{amsthm}
\usepackage{fancyhdr}
\usepackage{enumerate}
\usepackage[top=1in, bottom=1in, left=1in, right=1in]{geometry}
\usepackage{graphicx}
\usepackage[labelfont=bf]{caption}

\usepackage{caption}
\usepackage{subcaption}

\usepackage{verbatim}

\usepackage{wasysym}
\usepackage{hyperref}

\usepackage{wrapfig}

\begin{document}

\title{Frequency Spirals}
\author{Bertrand Ottino-L\"{o}ffler and Steven H. Strogatz}
\date{\today}

\begin{abstract}
We study the dynamics of coupled phase oscillators on a two-dimensional Kuramoto lattice with periodic boundary conditions. For coupling strengths just below the transition to global phase-locking, we find localized spatiotemporal patterns that we call ``frequency spirals.'' These patterns cannot be seen under time averaging; they become visible only when we examine the spatial variation of the oscillators' instantaneous frequencies, where they manifest themselves as two-armed rotating spirals. In the more familiar phase representation, they appear as wobbly periodic patterns surrounding a phase vortex. Unlike the stationary phase vortices seen in magnetic spin systems, or the rotating spiral waves seen in reaction-diffusion systems, frequency spirals librate:  the phases of the oscillators surrounding the central vortex move forward and then backward, executing a periodic motion with zero winding number. We construct the simplest frequency spiral and characterize its properties using analytical and numerical methods. Simulations show that frequency spirals in large lattices behave much like this simple prototype. 
\end{abstract}

\maketitle

\begin{quotation}
Spirals have long been objects of fascination in many parts of human culture, from art and architecture to science and mathematics. Within nonlinear dynamics, spirals have been studied in connection with such diverse phenomena as spiral density waves in galaxies, spiral waves of electrical activity in the heart and nervous system, growth spirals caused by screw dislocations in crystals, and spiral patterns of florets on the heads of sunflowers, to name just a few.  Here we consider a spiral pattern of a type that, as far as we know, has not been discussed previously. Instead of a  pattern of density or phase, it is a pattern of instantaneous frequencies: hence, a ``frequency spiral.'' We found it in a two-dimensional array of coupled oscillators whose natural frequencies were not quite identical, a regime in which the oscillators could no longer maintain global frequency entrainment and instead self-organized into one or more frequency spirals. We explore the properties of frequency spirals by analytical and geometrical methods and illustrate them by videotaped simulations of their dynamics in time and space. 
\end{quotation}

\section{Introduction}

\subsection{The Kuramoto model}

\noindent Many natural and technological systems can be viewed as enormous collections of self-sustained oscillators. Examples include networks of neurons and heart pacemaker cells, chorusing crickets, congregations of fireflies, arrays of Josephson junctions, oscillating chemical reactions, applauding audiences, and generators in the power grid~\cite{kuramoto84,strogatz00,pikovsky03}. These systems often display a remarkable capacity for mutual synchronization: if the oscillators interact strongly enough, they can spontaneously fall into a collective rhythm. Then a macroscopic fraction of the system runs at the same instantaneous frequency, despite the inevitable differences in the natural frequencies of the individual oscillators. 

The most widely studied mathematical model of mutual synchronization is the Kuramoto model~\cite{kuramoto84,strogatz00,pikovsky03, strogatz03, acebron05, dorfler14, rodrigues15}. It consists of a large collection of phase oscillators, coupled attractively to one another through a sine function of their pairwise phase differences. The oscillators themselves are diverse; their natural frequencies are randomly distributed across the population according to some prescribed probability distribution. In this way, the model embodies the competition between order and disorder mentioned above: the sinusoidal coupling tends to synchronize the oscillators, whereas the randomness in their natural frequencies tends to desynchronize them. 

The governing equations of the Kuramoto model are given by 
\begin{equation}
\dot \theta_i = \omega_i + K \sum_{j\in\mathcal{N}(i)} \sin(\theta_j - \theta_i),
\label{eqn:governing}
\end{equation}
where $\theta_i$ is the phase of oscillator $i$, $\omega_i$ is its natural frequency, $K \ge 0$ is the coupling strength, and the sum is over the neighbors $\mathcal{N}(i)$ of oscillator $i$. The natural frequencies are selected independently at random from a unimodal, symmetric distribution $g(\omega)$. By going into a suitable rotating frame, one can always assume without loss of generality that this $g(\omega)$ has zero mean. 

An undirected graph $G$ encodes the connection topology for \eqref{eqn:governing}; an edge between two oscillators signifies that they are coupled. For simplicity, Kuramoto~\cite{kuramoto84} began with the case where $G$ was a complete graph. In physical terms, this meant that the interactions were all-to-all, corresponding to an infinite-range or mean-field approximation. By using an ingenious self-consistency argument, Kuramoto showed that the system exhibited a phase transition to collective synchronization at $K_c = 2/[ \pi g(0)]$. As $K$ was increased above $K_c$, more and more oscillators near the middle of the frequency distribution $g(\omega)$ became \emph{phase locked} (meaning that $\dot \theta_i(t) - \dot \theta_j(t) \equiv 0$ for all those oscillators $i, j$, for all time $t$). Meanwhile the oscillators in the tails of $g(\omega)$ remained desynchronized from the locked group and from each other. 

The next wave of work relaxed the all-to-all assumption by allowing oscillators to be connected to their nearest neighbors on a one-dimensional chain or ring, a two-dimensional square grid, or a higher-dimensional cubic lattice~\cite{ermentrout85, sakaguchi87, sakaguchi88, strogatz88a, strogatz88b, daido88, ermentrout92, paullet94, ostborn02, ostborn04, hong05, hong07, ostborn09, lee10}. More recently, many researchers have explored the Kuramoto model on more complex topologies, allowing for nonlocal coupling~\cite{panaggio15} and small-world, scale-free, or other network architectures~\cite{dorfler14, rodrigues15}.

\subsection{Lattices}

In what follows, we revisit the classic problem of the Kuramoto model on a two-dimensional square lattice. This system was introduced nearly thirty years ago by Sakaguchi, Shinomoto, and Kuramoto~\cite{sakaguchi87}, but many things about it still remain mysterious. 

For example, we don't know for sure whether it exhibits macroscopic synchronization for any fixed, finite value of the coupling strength $K$.  More precisely, consider the model \eqref{eqn:governing} on an $M \times M$ square grid with nearest-neighbor coupling and periodic boundary conditions, and with $g(\omega)$ given by a standard unit normal. Fix $K$ and ask: For $N = M^2 \gg 1$, does the system display a frequency-locked cluster of size $O(N)$? By a \emph{frequency-locked cluster} we mean a connected subset of oscillators having the same long-term \emph{average frequency}
\begin{equation}
\tilde \omega_i = < \dot\theta_i> = \lim_{T\to\infty} \frac{\theta_i(t_0 + T ) - \theta_i(t_0)}{T}, 
\label{eqn:average_velocity}
\end{equation}
where $t_0$ is the transient time and $T$ is the averaging time \cite{sakaguchi87, sakaguchi88, strogatz88a, strogatz88b}. Numerical evidence \cite{sakaguchi87, hong05, hong07, lee10} and  renormalization arguments \cite{strogatz88a, strogatz88b,daido88,ostborn09} suggest that \emph{macroscopic} clusters (meaning clusters of size $O(N)$) cannot exist in the two-dimensional grid as $N \rightarrow \infty$ for any fixed $K$, but this has not been proven rigorously. The available evidence suggests that the coupling needs to grow like $K = O( \log N)$ for a macroscopic cluster to exist on the two-dimensional grid \cite{lee10}, at least for the values of $N$ that have been studied numerically. 

There are reasons to think that the two-dimensional case is marginal. In one-dimensional lattices, it has been proven that macroscopic clusters are impossible~\cite{strogatz88a, strogatz88b, daido88} for fixed $K$.  On the other hand, macroscopic clusters appear to exist in simulations~\cite{sakaguchi87, hong07} of lattices of dimension $d \ge 3$, as was suspected from the start on the basis of scaling and renormalization arguments \cite{sakaguchi87, strogatz88a, strogatz88b, daido88}.  But from a mathematical standpoint, the existence of macroscopic frequency-locked clusters for $d \ge 3$ remains conjectural. In fact, we don't even know whether the average frequencies ${\tilde \omega_i}$ are well-defined for Kuramoto lattices in any dimension; there is no proof that the long-time limits exist (although numerically, they certainly seem to). 

On the positive side, we do know that for any fixed value of $K$, if macroscopic clusters exist in dimension $d \ge 2$, they must be sponge-like at all length scales, in the sense that they cannot contain any macroscopic cubes of frequency-locked oscillators of size $ \epsilon N$, for any fixed $\epsilon > 0$, as $N \rightarrow \infty$ \cite{strogatz88a, strogatz88b}.  Finally, we also know that complete phase-locking is impossible in \emph{any} dimension, for fixed $K$, as $N$ tends to infinity \cite{strogatz88a, strogatz88b}.

\section{Phase vortices and frequency spirals}

Against this confusing backdrop, in 2010 Lee et al.~\cite{lee10} reexamined the two-dimensional Kuramoto lattice from a fresh perspective. They gave a neat (but non-rigorous) explanation for the observed $\log N$ scaling of the coupling $K_L$ at which the lattice becomes completely phase-locked. They thereby shed light on the same $\log N$ scaling seen for the coupling $K_E$ at the onset of entrainment (synonymous with what we have called macroscopic frequency-locking), and for the coupling $K_\text{chaos}$  at the onset of chaos and single phase slips.  They also illuminated the key role of what they called \emph{vortices} moving along the boundaries of the clusters in the system. To clarify these ideas, let us try to visualize the spatiotemporal state of a Kuramoto lattice.  

\subsection{Visualizing phases and average frequencies}

Following the seminal paper by Sakaguchi et al.~\cite{sakaguchi87}, the dynamics of Kuramoto lattices have traditionally been visualized in one of two ways: either one looks at the oscillators' phases $\left\{\theta_i(t) \right\}$ or their average frequencies $\left\{\tilde \omega_i \right\}$. 

\begin{figure}
%\begin{subfigure}[b]{0.50\textwidth}
%\includegraphics[width=\textwidth]{LowK_Phases.png}
%\caption{}
%\label{StandardMess(Phase)}
%\end{subfigure}
%~
%\begin{subfigure}[b]{0.50\textwidth}
%\includegraphics[width=\textwidth]{LowK_AveVel.png}
%\caption{}
%\label{StandardMess(AveVel)}
%\end{subfigure}

\includegraphics[width =0.3\textwidth]{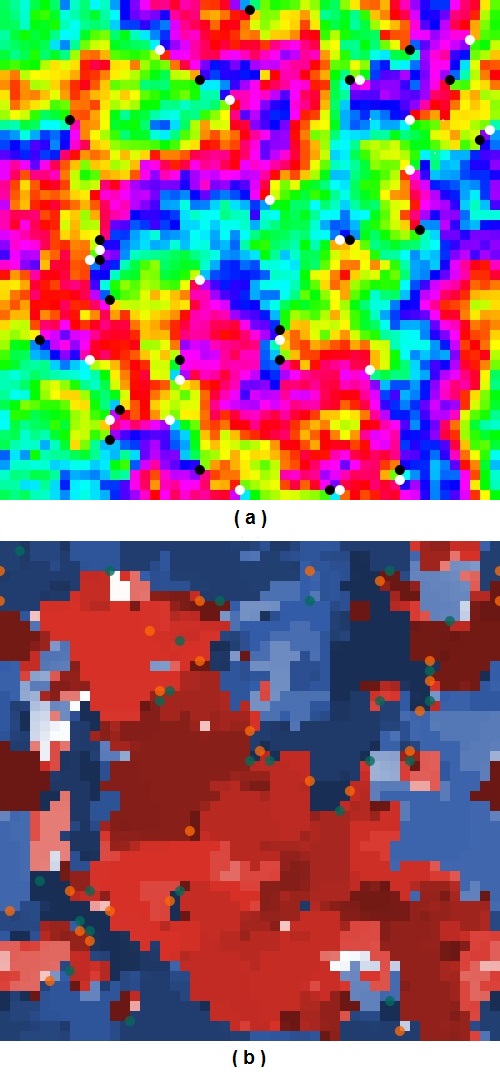}

\caption{Phase and average-frequency representations of a two-dimensional Kuramoto lattice, far below the transition to complete phase-locking ($K=1$, natural frequencies drawn from a standard normal distribution, uniform random initial phases, $50 \times 50$ grid with periodic boundary conditions). (a)~Phase representation: Colored pixels denote the instantaneous phases $\theta_i$, from 0 to 2$\pi$, of the oscillators according to a standard RGB color wheel. The pixel-sized white and black circular dots mark the location of phase vortices. \href{https://www.youtube.com/watch?v=OOYGXokE0wY}{\underline{Video 1}} shows the system's time-evolution in the phase representation (multimedia view).  (b)~Average-frequency representation: Color indicates the sign of the average frequency $\tilde \omega_i$, and darkness indicates its magnitude on a log scale. Blue pixels represent oscillators of positive average frequency, red pixels represent oscillators of negative average frequency, and white pixels represent oscillators with negligible average frequency (on the order of $10^{-4}$ in magnitude). The pixel-sized green and orange circular dots represent the sites of phase vortices, and are colored corresponding to their sign.  \href{https://www.youtube.com/watch?v=Xw8012hg0D8}{\underline{Video 2}} shows the evolution in the average-frequency representation  (multimedia view).  Notice how the disorder of the phase picture (a) gets smoothed over by the average-frequency picture (b). Figures and videos were made using Netlogo \cite{Wilensky}. }
\label{StandardMess}
\end{figure}

Figure~\ref{StandardMess}a (multimedia view) shows the phase representation for a $50 \times 50$ lattice far below the transition to complete phase-locking. It appears to be an inglorious mess. The same state looks much cleaner if we examine the average frequencies of the oscillators (Figure~\ref{StandardMess}b (multimedia view)). This way of looking at the lattice reveals large clusters of frequency-locked oscillators, shown as irregular blocks of a single color and darkness. 

%%%%%%  Figure 1 %%%%%%%%%

%%%%%%%%%%%%%%%%%%%%%%%%%%%%%%%%

%
%
%
%\begin{wraptable}{r}{2.0cm}
%\begin{tabular}{|c|c|}
%\hline
%0 & $\pi/4$ \\
%\hline
%$3\pi/2$ & $5\pi/6$ \\
%\hline
%\end{tabular}
%\caption{A demo phase vortex.}
%\label{PhaseVortex}
%\end{wraptable}

Although the average-frequency picture has its advantages, it misses some interesting features of the dynamics inherent in the phase picture. This is especially true if we look at a lattice that is almost completely locked to a single average frequency. Naturally, there are only a small number of features this close to total entrainment. This simplicity motivated Lee et al. \cite{lee10} to examine this case in detail. One particular structure that caught their attention was what they called a vortex, and what we will specifically refer to as a {\it phase vortex}. 

We say that a phase vortex exists on a $2 \times 2$ unit cell of adjacent oscillators if their phases exhibit a lattice curl of $\pm2\pi$ around the cell. The curl is easily computed by taking the sum of the ordered phase differences around the cell, after having shifted the phase differences mod $2\pi$ to be within $(-\pi, \pi]$. This sum can only ever be one of $-2\pi, +2\pi$ or $0$, which correspond to a negative phase vortex, a positive phase vortex, and no phase vortex, respectively.  We regard the phase vortex  as living at the center of the unit cell of four oscillators, rather than on any one oscillator. 

Phase vortices may be thought of as topological objects. They are conserved by the dynamics, and  vanish only if they collide with a vortex of the opposite sign or fall off the edge of a finite grid. So on a torus, we expect an equal number of positive and negative vortices, which can only ever be created or destroyed in pairs. 

Vortices have important implications for the average frequencies: they delineate the boundaries of frequency-locked clusters with their movement, as discovered by Lee et al.~\cite{lee10}. In one scenario, a single vortex circulates around the boundary of the cluster. In a second scenario, a pair of oppositely charged vortices is periodically created at a certain spot; the vortices then split apart, travel along opposite sides of the cluster, and finally annihilate each other at the far end of the cluster boundary.  In both cases, Lee et al.~\cite{lee10} found that when a phase vortex moves between two oscillators, the pair accumulates a $2\pi$ phase slip.  Such phase slips cause the differences in average frequency across the boundary between two clusters. 

The motion of the vortices also explains why the clusters are sometimes stable and sometimes unstable.  The stable case corresponds to the vortex-pair scenario described above. The consistent motion of the paired vortices leads to clusters with time-independent boundaries. In contrast, at lower coupling strength, the vortices move irregularly. As they wander around chaotically, they chop the macroscopic cluster into shreds. In this way, the onset of chaotic vortex motion goes hand in hand with the loss of frequency order in the system. 

%Additionally, a little algebra proves that at least one of the four phase differences in the unit cell of a vortex is at least $\pi/2$ in magnitude. Considering that we are using sine-based coupling, phase differences of larger than $\pi/2$ in magnitude correspond to a regime of negative coupling slope. Both these qualities make them portentous for a nearest-neighbor Kuramoto system. 

Finally, if we extend a phase vortex to the entire grid (thus forcing the grid to have a nonzero winding number around its boundary), we get a {\it phase spiral}. These spatially-extended structures have been proven to be stable phase-locked solutions under Kuramoto dynamics in the case where all the oscillators have identical frequencies~\cite{paullet94}.

\subsection{Visualizing instantaneous frequencies}

\begin{figure}
%\begin{subfigure}[b]{0.50\textwidth}
%\includegraphics[width=\textwidth]{SmallSpirals_Phases.png}
%\caption{}
%\label{SmallSpirals(Phase)}
%\end{subfigure}
%~
%\begin{subfigure}[b]{0.50\textwidth}
%\includegraphics[width=\textwidth]{SmallSpirals_AveVel.png}
%\caption{}
%\label{SmallSpirals(AveVel)}
%\end{subfigure}
%
%\centering
%\begin{subfigure}[b]{0.50\textwidth}
%\includegraphics[width=\textwidth]{SmallSpirals_Vel.png}
%\caption{}
%\label{SmallSpirals(Vel)}
%\end{subfigure}

\includegraphics[width =0.25\textwidth]{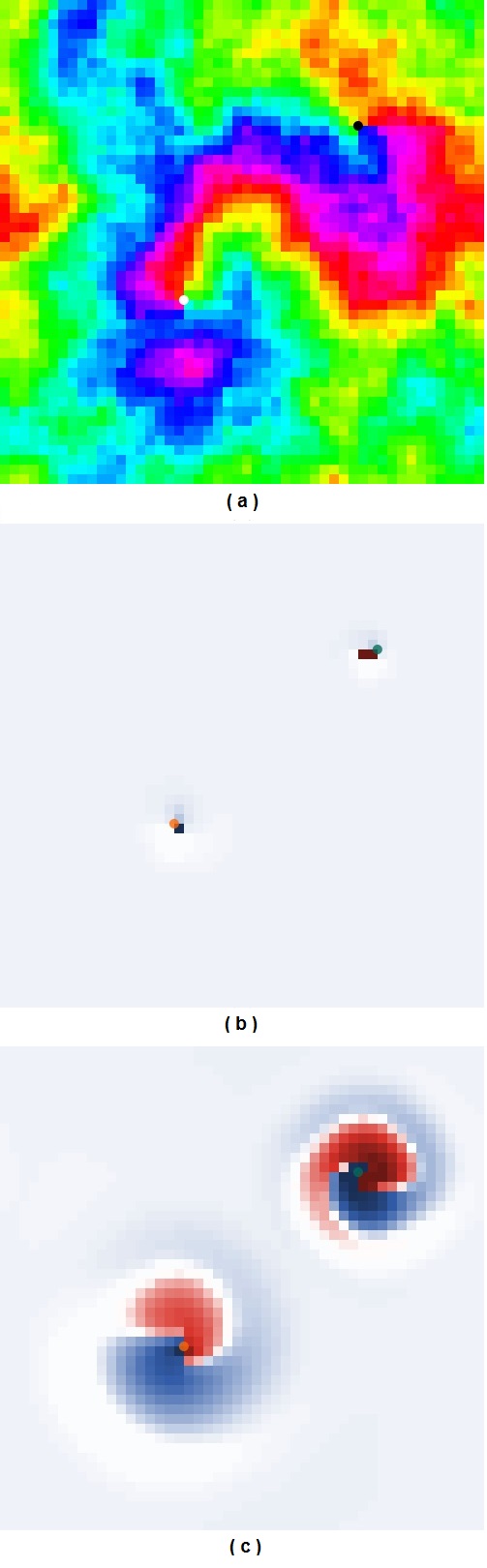}

\caption{Three different representations of a $50 \times 50$ Kuramoto lattice close to a phase-locked state (Parameters: $K=1.80$, natural frequencies drawn from a standard normal distribution, uniform random initial phases.) (a)~Phase representation, with phase $\theta_i$ from $0$ to $2\pi$ sent to an RGB color wheel.  \href{https://www.youtube.com/watch?v=I4mNXDrXVRE}{\underline{Video 3}} shows a movie of the lattice's time-evolution in the phase picture (multimedia view). (b)~Average-frequency representation, with color indicating sign of $\tilde \omega_i$ and darkness showing magnitude of $\tilde \omega_i$ on a log scale. \href{https://www.youtube.com/watch?v=xGfYK5v6yRA}{\underline{Video 4}} shows the evolution in the average-frequency picture after skipping a transient (multimedia view).   (c)~Instantaneous-frequency representation, with color indicating sign of $\dot \theta_i$ and darkness showing magnitude of $\dot \theta_i$ on a log scale. This representation reveals two frequency spirals. \href{https://www.youtube.com/watch?v=rlb0uDAPI88}{\underline{Video 5}} shows the evolution in the instantaneous-frequency picture (multimedia view).  In each panel, the pixel-sized circular dots denote the location of phase vortices, with colors representing the sign of the vortex. Notice how featureless the average-frequency picture (b) is compared to the phase picture (a) and the instantaneous-frequency picture (c). Figures and videos were made using Netlogo \cite{Wilensky}.}
\label{SmallSpirals}
\end{figure}

The visualization techniques discussed so far have ignored the oscillators' \emph{instantaneous frequencies} $\{\dot \theta_i(t)\}$, or equivalently, their  velocities as they move around their phase circles. Figure~\ref{SmallSpirals} shows why it can be illuminating to consider them. The simulation setup is similar to that of Fig.~\ref{StandardMess}, except the random initial conditions $\{\theta_i(0)\}$ and the random natural frequencies $\{\omega_i\}$ came from a different random seed, and the coupling was taken far closer to the onset of total phase-locking. Figure~\ref{SmallSpirals}a (multimedia view) and  \href{https://www.youtube.com/watch?v=I4mNXDrXVRE}{\underline{Video 3}} show that there is still a large amount of phase disorder, but the average-frequency picture (Fig.~\ref{SmallSpirals}b (multimedia view)) show that almost all the oscillators ultimately form a single macroscopic cluster, with only three oscillators still running at a different average frequency. What is happening in the neighborhood of these ``defecting" oscillators? 

%%%%%%% Figure 2  %%%%%%%%%%%

%%%%%%%%%%%%%%%%%%%%%%%%%%%%

As it turns out, interesting things are happening, as shown in Fig.~\ref{SmallSpirals}c (multimedia view). If you watch this system evolve from the perspective of its instantaneous frequency field, as in Fig.~\ref{SmallSpirals}c (multimedia view), you will notice that regions of extreme gradients in $\dot \theta$ tend to occur across phase vortices, with strong positive $\dot \theta$ on one side and strong negative $\dot \theta$ on the other. These dipoles move around the grid as bloopy little jellyfish-like objects, demarcating regions of differently changing phases. They survive the transient, causing the conspicuous spiral structures in  Fig.~\ref{SmallSpirals}c (multimedia view). We can clearly identify two spinning arms coming out of a central core: one of positive $\dot \theta$, the other negative. Loosely speaking, it nearly appears as if we have a pair of rotating dipoles. Because they rotate at a constant rate and the two arms are equally intense, these arms end up averaging out in the average-frequency picture (Fig.~\ref{SmallSpirals}b (multimedia view)). And that is probably why these structures  have not been reported before. Under time averaging, they become invisible for the same reason that $\sin t$ would be invisible: they both  average to zero. 

We call these structures {\it frequency spirals}, since they appear as spirals in the instantaneous frequency field. Lee et al.~\cite{lee10} have implicitly argued that frequency spirals are important to understanding the entrainment threshold of the nearest-neighbor Kuramoto model, so let us devote some time to understanding their phenomenology.

\section{The canonical frequency spiral}

A frequency spiral on a $50 \times 50$ lattice involves $5 \times 10^3$ tunable parameters: $50 \times 50$ initial phases and $50 \times 50$ natural frequencies. So one might worry that it takes a bit of contrivance to make a frequency spiral happen. 

\begin{figure}[h]
%\begin{subfigure}[b]{0.50\textwidth}
%\includegraphics[width=\textwidth]{CanonSpiral_Phases.png}
%\caption{}
%\label{CanonlSpirals(Phase)}
%\end{subfigure}
%~
%\begin{subfigure}[b]{0.50\textwidth}
%\includegraphics[width=\textwidth]{CanonSpiral_Vel.png}
%\caption{}
%\label{CanonSpirals(AveVel)}
%\end{subfigure}

\includegraphics[width =0.3\textwidth]{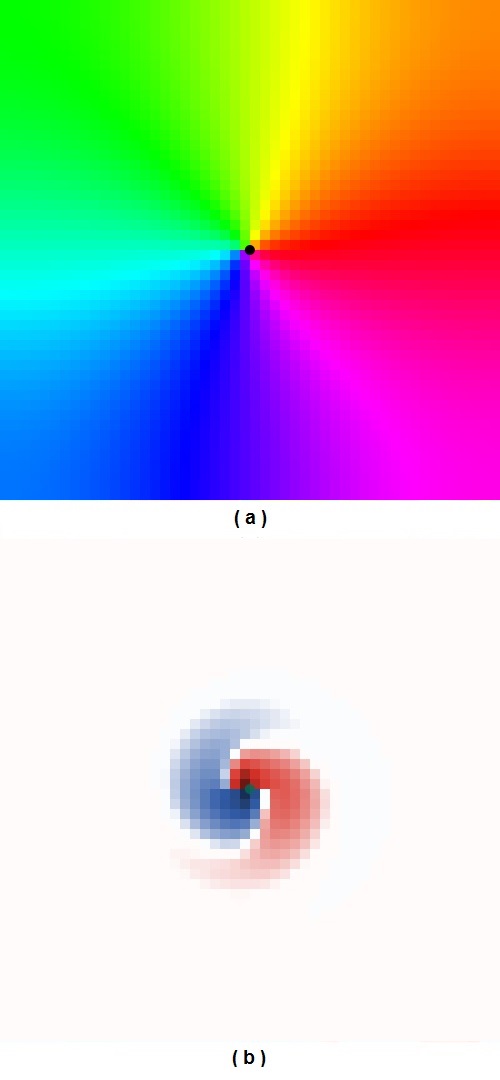}

\caption{The canonical frequency spiral, shown in (a) phase and (b) instantaneous-frequency representations. Parameters for numerical integration of Eq.~\eqref{eqn:governing}: coupling $K=1.0$; all natural frequencies set to zero, except for a central one with $\omega = 0.75$; initial phases chosen to form a phase spiral. Notice how it is impossible to discern the two arms of the frequency spiral in the phase picture (a). To watch the system's time-evolution in the phase picture or instantaneous-frequency picture, click \href{https://www.youtube.com/watch?v=naE_v-mRbgw}{\underline{Video 6}} (multimedia view). Figures and videos were made using Netlogo \cite{Wilensky}.}
\label{CanonSpirals}
\end{figure}

Thankfully, this isn't the case! The creation of a frequency spiral is nearly as easy as it looks. To construct the simplest possible example -- the \emph{canonical frequency spiral} shown in Fig.~\ref{CanonSpirals} (multimedia view) -- just place a phase vortex next to a single oscillator of sufficiently large natural frequency. Set all other oscillators to have $\omega_i = 0$, and normalize the coupling in Eq.~\eqref{eqn:governing} to $K = 1$ by rescaling time. Let $\omega$ denote the frequency of the distinguished central oscillator in these time units. To make things even simpler, we can take the topologically-required second phase vortex and shove it to the point at infinity. Doing this means that we can no longer work on a torus, since the phase pattern must have a winding number of 1 around the boundary of the grid. So we will use free boundaries from now on, unless otherwise noted. Choose the initial phases $\{\theta(x, y , t=0)\}$ so that they have a phase vortex at the center of the lattice; a natural choice satisfies $\cos \theta(x, y , t=0)  = x$ and $ \sin \theta(x, y , t=0)  = y.$ The resulting pattern of natural frequencies and initial phases is reminiscent of the phase spiral constructed by Paullet and Ermentrout~\cite{paullet94}, except one of the central four oscillators now has a natural frequency different from all the others.

%%%%%       Figure 3      %%%%%%%%%%

\subsection{Spatial decay of the frequency spiral}

Having constructed a canonical frequency spiral, it is time to numerically drill down into its features. Numerical integration of \eqref{eqn:governing} shows that rotation period of the spiral doesn't change much from cycle to cycle. In this sense, the structure is persistent in time. 

\begin{figure}[h]
%\begin{subfigure}[b]{0.50\textwidth}
%\includegraphics[width=\textwidth]{SpiralAmpDecay.pdf}
%\caption{} 
%\label{SpiralAmpDecay}
%\end{subfigure}
%~
%\begin{subfigure}[b]{0.50\textwidth}
%\includegraphics[width=\textwidth]{CanonSpiralBif.png}
%\caption{}
%\label{SpiralFreqBif}
%\end{subfigure}

\includegraphics[width =0.45\textwidth]{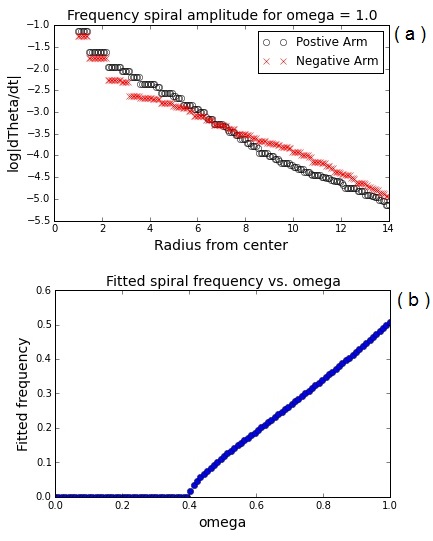}

\caption{Properties of the canonical frequency spiral. (a) Spatial decay of the spiral arms. The amplitude of the instantaneous frequency decays roughly exponentially with distance from the center. (b) Rotation rate of the frequency spiral versus magnitude of $\omega$ in the center. Plots were made using Numpy and Matplotlib \cite{Numpy, Matplotlib}. The computations were performed on grids of $149 \time 149$ oscillators, using a fourth-order Runge Kutta method with a timestep of 0.125. Figure \ref{SpiralGraphs}a was made using a transient time of 32000 time units. Figure \ref{SpiralGraphs}b was made using a transient period of 4000, and the period was measured over a length of 24000 time units.}
\label{SpiralGraphs}
\end{figure}

In space, however, it decays. We can take a freeze-frame and examine how its positive and negative arms of high and low instantaneous frequency attenuate with distance from the core (Fig.~\ref{SpiralGraphs}a). As suggested by prior figures, both arms die off nearly exponentially in space, and are of equal strengths. However, the semi-log plot doesn't exactly look like two overlapping straight lines. This indicates that there are some complications, likely from the discreteness of the grid, the finiteness of the simulation, and most importantly, the nonlinearity of sine. In any case, this exponential-like decay goes a long way to explaining why these structures were undetectable in the phase view: most of the action was too small to see! If we look far from the origin and take a continuum limit, the small phase differences mean that our dynamics start to resemble the heat equation with a diffusion constant of exactly $K$. So this sort of exponential decay is to be expected.

%%%%%%%    Figure 4    %%%%%%%%%%%%%%%

%%%%%%%%%%%%%%%%%%%%%%%%%%%%%

\subsection{Rotation rate of the frequency spiral}

Since we normalized $K=1$ and preselected all the initial phases and natural frequencies, the canonical frequency spiral has exactly one tunable parameter: $\omega$. It alone determines the rotation rate of the spiral. 

The resulting functional dependence, computed from simulations, is plotted in Figure~\ref{SpiralGraphs}b. For large values of $\omega$, the relationship is nearly linear, as one might expect. At sufficiently small $\omega$, the frequency spiral doesn't exist at all, and we just get a static phase spiral. This is also to be expected, since it has been proven that when all the natural frequencies are identically zero, static phase spirals are linearly stable solutions~\cite{paullet94}. Interpreting these solutions as hyperbolic fixed points of our dynamical system, one expects them to persist for sufficiently small perturbing natural frequencies. 

Putting these observations together gives us the overall shape of Figure~\ref{SpiralGraphs}b: a zero branch with a curve bifurcating from it. The bifurcation point, $\omega = \omega_c \approx 0.40$, therefore corresponds to a frequency spiral with zero rotation rate, or equivalently, infinite rotation period. Recalling that $\omega$ actually represents the difference between the central oscillator's natural frequency and those of the oscillators around it, we infer (as a rule of thumb) that in a large, two-dimensional Kuramoto lattice, a phase spiral will become a frequency spiral if it gets caught on an oscillator whose natural frequency is $0.40$ larger than that of its neighbors. 

\begin{figure}
\includegraphics[width =0.41\textwidth]{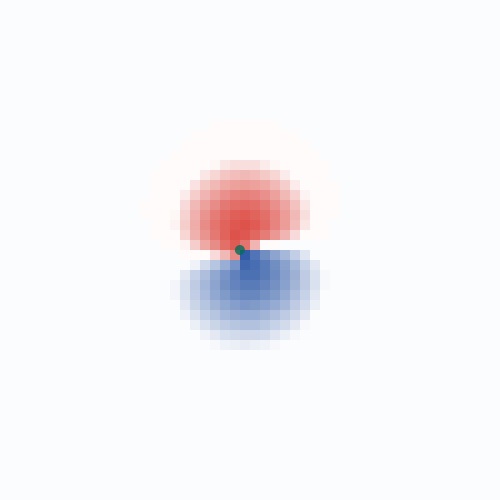}
\caption{A still of \href{https://youtu.be/WF7PqrPwXcU?t=30s}{\underline{Video 7}} (multimedia view), which shows the dynamics of the canonical frequency spiral at $\omega$ = 0.41, close to the bifurcation. In a single cycle, the spiral rotation slows down four times, corresponding to the bottlenecks left behind by the four SNIPERs.  Figure and video were made using Netlogo \cite{Wilensky}.}
\label{NearCrit}
\end{figure}

Because of this expected infinite period, it seems plausible that the birth of the frequency spiral involves a Saddle-Node Infinite PERiod bifurcation (also known as a SNIPER). If we look closely at the simulated dynamics of the frequency spiral near the critical $\omega$, as in Figure~\ref{NearCrit} and \href{https://youtu.be/WF7PqrPwXcU?t=30s}{\underline{Video 7}} (multimedia view), we see the spiral nearly stops \emph{four} times in a single period, once for each grid direction. In fact, it stops looking very spiral-like at all, since the pauses cause the arms themselves to mostly decay away as the system relaxes. This leads to complications when studying spiral characteristics (e.g., decay, twist) near the bifurcation. But these four pauses hint that we are seeing four slow regions in state space, each corresponding to the ghost of a colliding saddle and node. As such, we actually expect there to be {\it four} SNIPERs occurring simultaneously! However, the four-way symmetry may be slightly broken, given that we chose only one of the four oscillators in the center of the phase vortex to carry the $\omega$, imparting a slight offset to the $\pi/2$ symmetry natural to the phase spiral.

%%%%%%%%%%%%%% Figure 5    %%%%%%%%%%%%%%

%%%%%%%%%%%%%%%%%%%%%%%%%%%%%%%%%%

%\newpage

\section{Minimal model of a frequency spiral}

The numerical evidence above suggests that the canonical frequency spiral is created by four simultaneous SNIPER bifurcations. To confirm this scenario, we would need to analyze the existence and stability of a static phase spiral as we vary $\omega$. Unfortunately no closed-form solution for a static phase spiral on a two-dimensional Kuramoto lattice is known, except for the case $\omega = 0$~\cite{paullet94}.

To make progress, suppose we simplify the frequency spiral even further. As it turns out, we can get away with surprisingly few details.  Five active oscillators are enough to make a recognizable frequency spiral. (The fewer the oscillators the better, considering that even modestly sized grids of $10 \times10$ oscillators can have upwards of $10^5$ fixed points \cite{mehta15}.) 

\subsection{Constructing the minimal model}
\begin{figure}
%\begin{subfigure}[b]{0.50\textwidth}
%\includegraphics[width=\textwidth]{FiveTile_Phases.png}
%\caption{}
%\label{5Tile(Phase)}
%\end{subfigure}
%~
%\begin{subfigure}[b]{0.50\textwidth}
%\includegraphics[width=\textwidth]{FiveTile_Vel.png}
%\caption{}
%\label{5Tile(Vel)}
%\end{subfigure}

\includegraphics[width =0.3\textwidth]{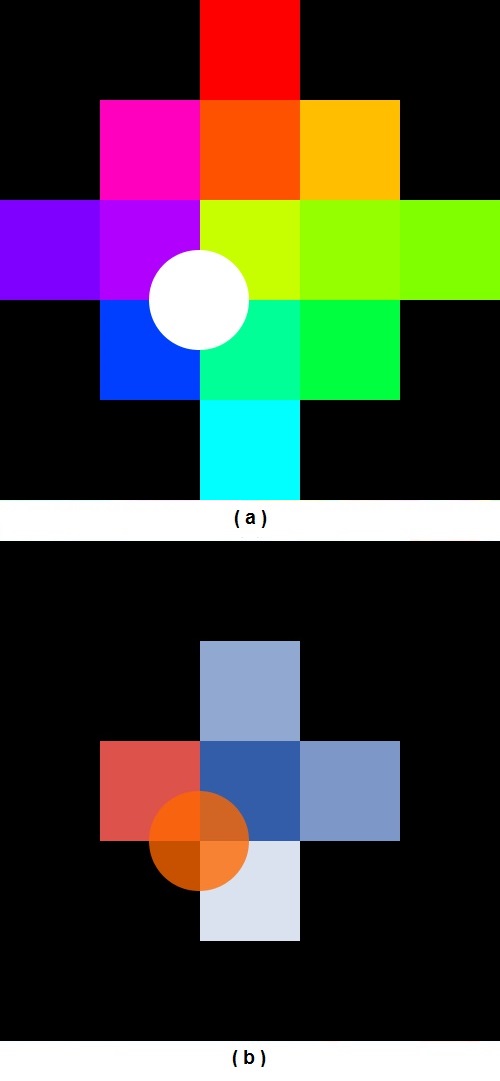}

\caption{Visualization of the minimal model~\eqref{eqn:5TileEqs} in (a) phase and (b) instantaneous-frequency representations. Only the central five oscillators are active. The other oscillators are held at fixed phases, and therefore appear black in the instantaneous-frequency picture (b). The circular dot in (a) and (b) marks a phase vortex. Click \href{https://www.youtube.com/watch?v=FCjJfuYO1_A}{\underline{Video 8}} to watch the system's time evolution (multimedia view). Figures and videos were made using Netlogo \cite{Wilensky}.}
\label{5Tile}
\end{figure}

Recall the recipe to make the canonical frequency spiral:  put a high $\omega$ oscillator in the middle of a phase spiral. The minimal way to do this is shown in Table~\ref{5TileTable} and Figure~\ref{5Tile} (multimedia view). 

%%%%%%%%%%%%%%      Table 1       %%%%%%%%
\begin{table}[htb]
\begin{center}
\begin{tabular}{|c|c|c|c|c|}
\hline
& & $ 0 $ & & \\
\hline
& $\pi/4 $ & $\theta_0 + 0 $ & $7\pi/4$ & \\
\hline
$\pi/2$ & $\theta_1 + \pi/2$ & $\zeta $ & $\theta_3 + 3\pi/2 $ &$3\pi/2$ \\
\hline
& $3\pi/4$ & $\theta_2 + \pi $ & $5\pi/4 $ & \\
\hline
& & $\pi$ & & \\
\hline
\end{tabular}
\end{center}
\caption{Definition of the oscillators' phases in the minimal model. Notice that the variables $\theta_i$ have been redefined here; they now denote deviations from the constant phase offset found in a static phase spiral.}
\label{5TileTable}
\end{table}
%%%%%%%%%%%%%%%%%%%%%%%%%%%%%%%%%%

%%%%%%%     Figure 6       %%%%%%%%%%

%%%%%%%%%%%%%%%%%%%%%%%%%%%%

\noindent Here $\theta_i(t)$ and $\zeta(t)$ are dynamical variables, and are coupled to their four nearest neighbors. Their dynamics are inherited from  the Kuramoto model~\eqref{eqn:governing}. The surrounding oscillators are held at fixed phases (the dynamics do not apply to them). Fixing these outer oscillators is reasonable because it approximates what we know about the canonical frequency spiral:  the instantaneous velocities of the oscillators surrounding the core die off exponentially fast in space, so a completely motionless far-field is not too qualitatively different. This approximation also allows us to enforce the nonzero winding number without needing to use large grids. We use a special symbol $\zeta$ for the central oscillator, because it obeys a noticeably different equation from the $\theta_i$: 
\begin{alignat}{2} \label{eqn:5TileEqs}
& \dot\theta_0 =&& \sin(\zeta - \theta_0) -\kappa \sin\theta_0, \notag \\
& \dot\theta_1 =&& \sin(\zeta - \theta_1 - \pi/2) -\kappa \sin\theta_1, \notag \\
& \dot\theta_2 =&& \sin(\zeta - \theta_2 - \pi) -\kappa \sin\theta_2, \\
& \dot\theta_3 =&& \sin(\zeta - \theta_3 - 3\pi/2) -\kappa \sin\theta_3, \notag \\
& \dot\zeta =&& \omega - \sin(\zeta - \theta_0) - \sin(\zeta - \theta_1 - \pi/2) \notag \\
 & \hspace{0.5cm} &&- \sin(\zeta - \theta_2 - \pi) - \sin(\zeta - \theta_3 - 3\pi/2), \notag
\end{alignat}
where the coefficient $\kappa = 1+\sqrt{2}$ came out from some trigonometric identities to produce the rather symmetric dynamics seen above. For example, notice that the dynamics remain identical under rotations of $\zeta \to \zeta + \pi/2$ and $\theta_i\to \theta_{i+1}$. Ergo, if one equilibrium undergoes a bifurcation, then three other equilibria, related to it by this rotational symmetry, do so at the same time. So the simultaneity of the four  SNIPER bifurcations mentioned above has a nice explanation in this setting.

\subsection{Analyzing the minimal model}
\begin{figure}
%\begin{subfigure}[b]{0.50\textwidth}
%\includegraphics[width=\textwidth]{5TileBifDiagram_Small.pdf}
%\caption{}
%\label{5TileBif}
%\end{subfigure}
%~
%\begin{subfigure}[b]{0.50\textwidth}
%\includegraphics[width=\textwidth]{5TileContinutation.png}
%\caption{}
%\label{5TileContinuation}
%\end{subfigure}

\includegraphics[width =0.45\textwidth]{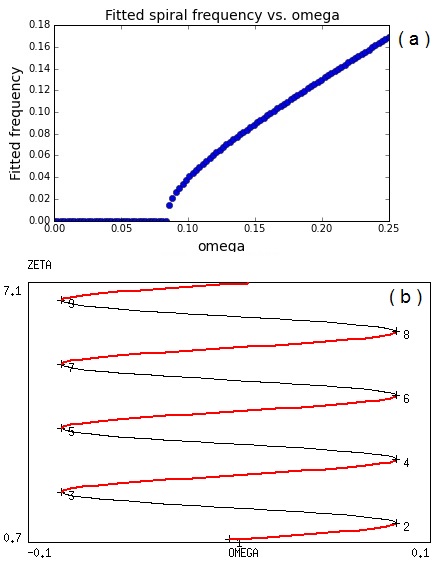}
\label{5TileGraphs}

\caption{Properties of the minimal model~\eqref{eqn:5TileEqs}. (a)~Rotation rate of a frequency spiral, plotted versus the natural frequency $\omega$ of the central oscillator. Plots were made using Numpy and Matplotlib \cite{Numpy, Matplotlib}. Simulations were generated using a timestep of 0.125, with a transient time of 500. Spiral rotation rates were measured over a span of 2000 time units. (b)~The equilibrium value of $\zeta$ is plotted versus $\omega.$ The red curves represent stable equilibria, and the black curves represent unstable equilibria. Four simultaneous SNIPER bifurcations occur at $\omega_c \approx \pm 0.08.$ Plot generated using XPP/AUTO, with a maximum stepsize of 0.05 and a minimum of 0.001.}
\label{5TilePlots}
\end{figure}

The minimal model~\eqref{eqn:5TileEqs} has a stable equilibrium\footnote{This is the only equilibrium with $\theta_0 = -\theta_1$ and $\theta_2 = -\theta_3$, up to some rotations. Fun fact: the solutions that have $\theta_0 = 0$ are actually unstable.} at $\omega = 0$ at $\zeta = \pi/4$, $\theta_0 = \arctan[1/(1+\sqrt{2}\kappa)] = -\theta_1$, and $\theta_2 = \arctan[1/(1-\sqrt{2}\kappa)]=-\theta_3$. As it turns out, for sufficiently large $\omega$, this equilibrium bifurcates and gives rise to the world's tiniest frequency spiral, as seen in Figure~\ref{5Tile}b and \href{https://www.youtube.com/watch?v=FCjJfuYO1_A}{\underline{Video 8}} (multimedia view). Since each arm is no more than one oscillator long, discussing the spatial decay of its amplitude is out of the question. But we can once again plot the rotation rate of the frequency spiral versus $\omega$ to get the familiar-looking plot in Figure~$\ref{5TilePlots}a$, and find a critical $\omega$ of about 0.08 for the minimal model. The reason this critical $\omega$ is so much smaller than in the full grid (0.08 vs. 0.40) is clear physically: by holding all the surrounding oscillators fixed, we made the system terribly stiff and brittle. As such, it is far more likely to break than bend, hence the low threshold for bifurcation.

%%%%%%%     Figure 7       %%%%%%%%%%%%

%%%%%%%%%%%%%%%%%%%%%%%%%%%

To determine the nature of the bifurcation, we can plug the five differential equations of the minimal model into our favorite continuation method. The bifurcation diagram shown in Figure \ref{5TilePlots}b was generated using XPP by starting off with the previously identified fixed point and continuing for a length of $2\pi$ in the $\zeta$ variable. We have plotted $\zeta$ (as opposed to any of the $\theta_i$) because it is the most visually interesting choice. We can also generate the same diagram by solving for the relevant equilibrium state of~\eqref{eqn:5TileEqs} and applying trig identities to get the following formula relating $\zeta$ to $\omega$: 
\begin{alignat}
 \omega \sqrt{\frac{2}{\kappa}} = &  \frac{\sin\zeta}{(\sqrt{2}+\cos\zeta)^{1/2}} + \frac{-\cos\zeta}{(\sqrt{2}+\sin\zeta)^{1/2}} \notag \\ 
&+\frac{-\sin\zeta}{(\sqrt{2}-\cos\zeta)^{1/2}} + \frac{\cos\zeta}{(\sqrt{2}-\sin\zeta)^{1/2}} . 
 \label{eqn:zeta_vs_omega}
 \end{alignat}
\noindent
As we expected, we get four simultaneous saddle-node bifurcations on the positive $\omega$ side and four on the negative side. These occur at $\omega_c \approx \pm 0.0833$. For $\omega \gg \omega_c$, the central oscillator $\zeta$ just heedlessly accumulates phase. For $\omega$ close to but above $\omega_c$, $\zeta$ still accumulates phase, but moves slowly at positions corresponding to the limit points of the bifurcation diagram. And for lower $\omega$, $\zeta$  gets stuck and the entire system becomes static.

%%%%%%%%%%%%%%%%%%     Figure  8    %%%%%%%%%%%%%%%%%%
\begin{figure}
\includegraphics[width =0.41\textwidth]{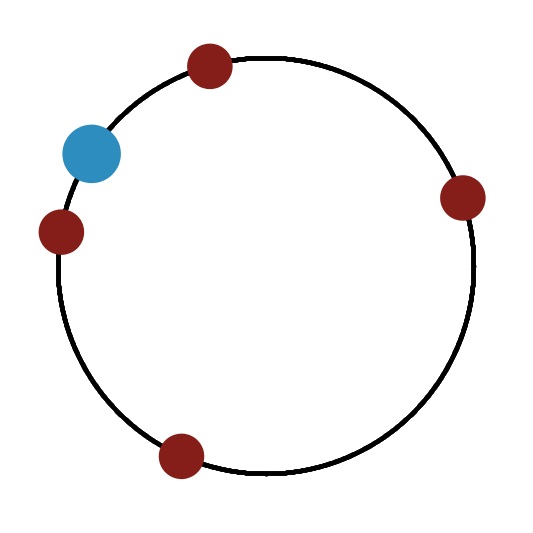}
\caption{The dynamics of  the minimal model in a state space representation. The image is a still frame from \href{https://www.youtube.com/watch?v=BNfUZkSQVJc}{\underline{Video 9}} (multimedia view), which shows the evolution of the oscillators' phases for $\omega$ = 0.25, in a regime just past the SNIPER bifurcations.  The large blue circle shows the circular phase of the central $\zeta$ oscillator, and the smaller red dots show the phases of its four neighbors.  Figure and video were made using Netlogo \cite{Wilensky}.}
\label{5TileCirc}
\end{figure}

%%%%%%%%%%%%%%%%%%%%%%%%%%%%%%%%%%%%%%%%%%%%

Figure~\ref{5TileCirc} and \href{https://www.youtube.com/watch?v=BNfUZkSQVJc}{\underline{Video 9}} (multimedia view) show an extremely revealing state-space view of what the dynamics are like after the SNIPERs have occurred and created a frequency spiral. Notice how the phase $\zeta$ runs monotonically around the circle, while the phases of the neighboring oscillators librate back and forth about their former equilibrium positions (which no longer exist, but which would have been stable in a static phase spiral).  This back and forth jiggling accounts for the  regions of positive and negative instantaneous frequency that make up the two arms of a frequency spiral.

That is really all that is going on in a frequency spiral: a high-frequency central oscillator breaks loose when locking is lost. Then, as it sweeps past the phases of its neighboring oscillators, it periodically tugs them forward and backward, thanks to the sinusoidal coupling. This forward and backward motion gives rise, in the instantaneous-frequency representation, to the two spinning arms of the spiral. Although this intuitive picture provides a highly simplified caricature of a frequency spiral in the setting of a minimal model, its fundamental mechanics should apply all the same to the full model too. 

\subsection{Perturbative results}
 
Even for the minimal model, the nonlinearities in the dynamics make exact results hard to come by. A few perturbative results of interest can be obtained for $\omega \gg \omega_c$. In that limit we can use Lindstedt's method~\cite{drazin} to approximate the periodic solution originally produced by the SNIPER bifurcations, and thereby obtain an estimate of the rotation rate of the frequency spiral in the regime far above threshold. 

To apply Lindstedt's method, we rescale time to a new variable $\tau = \Omega \epsilon t$, where $\epsilon := 1/\omega \ll 1$. This introduces a factor of $\Omega$ in front of our time derivatives. Then we seek solutions that are $2\pi$-periodic with respect to $\tau$, eliminating secular terms as needed to enforce this periodicity condition,  and in this way obtain the angular frequency $\Omega$ of the periodic solution, as well as solutions for $\theta_i$ and $\zeta$ that have $2\pi$-periodic oscillatory terms in the $\tau$ variable.  The method involves expanding $\Omega$, $\zeta$, and $\theta_0$ in terms of $\epsilon$, then equating the coefficients of like powers in $\epsilon$. We also choose to shift the origin of time so that $\zeta(0) = 0$, for the sake of clarity. 

Recalling that $\kappa = 1 + \sqrt{2}$, and omitting the analytical details since Lindstedt's method is standard~\cite{drazin}, we find the following asymptotic series: 
\begin{alignat}{1} \label{5TilePerturb}
\Omega &= 1 - 2 \epsilon^2 +  \epsilon^4 \left( 2\kappa^2 - \frac{7}{2} \right) + \mathbb{O}(\epsilon^5), \notag \\
\zeta(\tau) &= \tau + \epsilon^4 \left( \frac{1}{8}\sin4\tau \right)  + \mathbb{O}(\epsilon^5), \\
\theta_0(\tau) &= \epsilon \left(\frac{1}{2 \kappa} - \cos\tau \right) + \mathbb{O}(\epsilon^2) \notag
  % + \epsilon^2 \left(\kappa -\frac{1}{2\kappa}  +\frac{1}{2}\cos\tau \right) \sin\tau  + \mathbb{O}(\epsilon^3), \notag
\end{alignat}

\noindent and naturally, by the rotational symmetry, $\theta_1(\tau) = \theta_0(\tau-\pi/2), \theta_2(\tau) = \theta_0(\tau-\pi), \theta_3(\tau) = \theta_0(\tau+\pi/2)$. 

The first things to notice here are that $\zeta$ is linear in time, to lowest order; this corresponds to the monotonic running of the central oscillator around the phase circle seen in \href{https://www.youtube.com/watch?v=BNfUZkSQVJc}{\underline{Video 9}}. Also, notice that $\theta_0(\tau)$ is a small back-and-forth oscillation. This is the libration, the jiggling motion of the neighboring oscillators, also seen in the video.

It is interesting to note how small the corrections to $\zeta$ are, compared to the corrections in the other two quantities above. The deviations from constant linear growth in $\zeta$ only occur at $\epsilon^4$, compared to the $\epsilon$-sized amplitude of the $\theta_i$'s or the $\epsilon^2$ deviation of the frequency. Given the straightness of the bifurcating branch in Figure~\ref{5TilePlots}a, this near-linearity was to be expected. Moreover, the time-shift symmetry and the $\sin 4\tau$ dependence could have been guessed from the original symmetries in the equations \eqref{eqn:5TileEqs}. That is to say, if all the $\theta_i$'s obeyed some function of the form $\phi( \tau- i \pi / 2)$, then we find that we want $\zeta(\tau) = \tau + \Delta(4\tau)$, where $\Delta$ is $2\pi$-periodic. Since the rotation rate of a frequency spiral is determined by the motion of $\zeta$, this is all to say that we typically expect its rate to scale almost exclusively with $\omega$, and for the phases of nearby oscillators to vary as sinusoids (to leading order). These perturbative results are easy to check by simulation; all we need to do is numerically measure the frequency $\Omega$ as well the amplitudes of the oscillations in $\theta_0$ and $\zeta$. 

\begin{figure}
%\centering
%\begin{subfigure}[b]{0.75\textwidth}
%\includegraphics[width=\textwidth]{OmegaZetaPlotPreidct5Tile.pdf}
%\caption{The first correction in $\Omega$, as measured from the slope of $\zeta$. }
%\label{5Tile(OmegaZeta)}
%\end{subfigure}
%
%\begin{subfigure}[b]{0.75\textwidth}
%\includegraphics[width=\textwidth]{ThetaAmpPlotPredict5Tile.pdf}
%\caption{The amplitude of the oscillations in $\theta_0$. }
%\label{5Tile(ThetaAmp)}
%\end{subfigure}
%
%\begin{subfigure}[b]{0.75\textwidth}
%\includegraphics[width=\textwidth]{ZetaAmpPlotPredict5Tile.pdf}
%\caption{The amplitude of the oscillations in $\zeta$. }
%\label{5Tile(ZetaAmp)}
%\end{subfigure}

\includegraphics[width =0.45\textwidth]{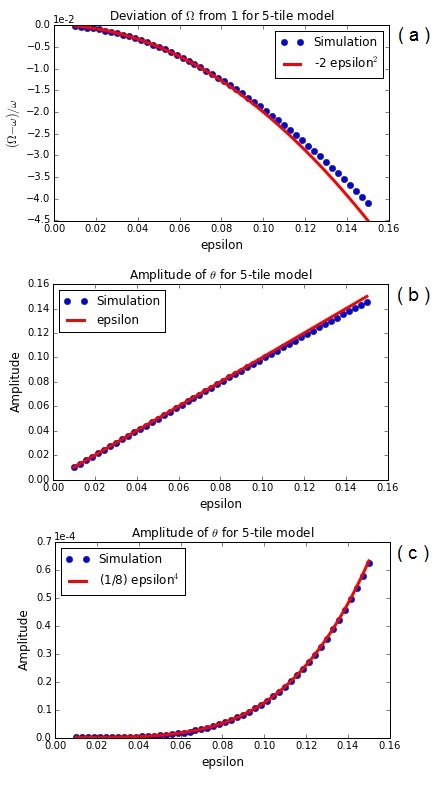}
\label{5TileLinst}

\caption{Comparison of perturbation theory and simulation of a frequency spiral for the minimal model of five oscillators given by Eq.~\eqref{eqn:5TileEqs}. Simulation results are plotted on top of the lowest-order prediction obtained by Lindstedt's method. At small $\epsilon$ the overlap is nice, as predicted. The simulations used timesteps of $\epsilon/50$, and used a transient time of about 500 periods, followed by an observation time of about 100 periods. Figures made using Numpy and Matplotlib \cite{Numpy, Matplotlib}.}
\label{5TileAsymptotics}
\end{figure}

Figure~\ref{5TileAsymptotics} shows that our predictions hold up well against simulation of the minimal model at small $\epsilon$, though higher-order corrections start to dominate for larger $\epsilon$.

%%%%%%%%%%%%%    Figure 9    %%%%%%%%%%%%%%%%
%%%%%%%%%%%%%%%%%%%%%%%%%%%%%%%%%%%%

However, we can do one better. Returning to the full grid, we can guess that the distinguished central oscillator will obey a similar behavior to $\zeta$, and a neighboring oscillator will behave like $\theta_0$. 

\begin{figure}
\centering
%
%\begin{subfigure}[b]{0.75\textwidth}
%\includegraphics[width=\textwidth]{OmegaZetaPlotPreidctFull.pdf}
%\caption{The first correction in $\Omega$, as measured from the slope of $\zeta$. }
%\label{Full(OmegaZeta)}
%\end{subfigure}
%
%\begin{subfigure}[b]{0.75\textwidth}
%\includegraphics[width=\textwidth]{ThetaAmpPlotPredictFull.pdf}
%\caption{The amplitude of the oscillations in $\theta$. }
%\label{Full(ThetaAmp)}
%\end{subfigure}
%
%\begin{subfigure}[b]{0.75\textwidth}
%\includegraphics[width=\textwidth]{ZetaAmpPlotPredictFull.pdf}
%\caption{The amplitude of the oscillations in $\zeta$. }
%\label{Full(ZetaAmp)}
%\end{subfigure}

\includegraphics[width =0.45\textwidth]{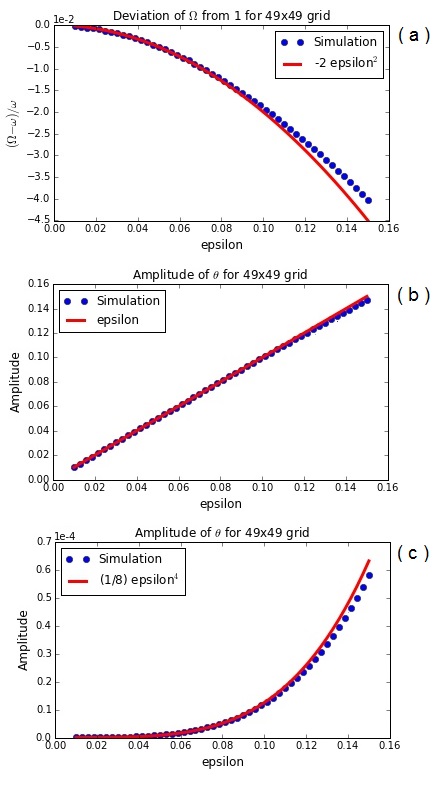}

\caption{Properties of a frequency spiral obtained via simulation of a $49 \times 49$ grid, compared to the lowest-order perturbative results obtained by Lindstedt's method applied to the five-oscillator model~\eqref{eqn:5TileEqs}. At small $\epsilon$ the overlap is good, despite this being an uncontrolled approximation. The simulations used timesteps of $\epsilon/50$, and used a transient time of about 500 periods, followed by an observation time of about 100 periods. Figures made using Numpy and Matplotlib \cite{Numpy, Matplotlib}.}
\label{FullAsymptotics}
\end{figure}

Figure~\ref{FullAsymptotics} shows these guesses are well founded, with the results from the Lindstedt analysis continuing to give a good approximation for small $\epsilon$, despite the vastly larger number of oscillators in the full grid. This suggests that our minimal model has qualitatively emulated the dynamics of frequency spirals in a much larger Kuramoto lattice.

%%%%%%%%    Figure 10       %%%%%%%%%%%%

%%%%%%%%%%%%%%%%%%%%%%%%%%%%%%%

\section{Discussion and future directions}

We have shown here that there is something to be learned by looking at Kuramoto lattices in a new way, using instantaneous frequencies as a window into their dynamics. It is understandable why no one did this in the distant past, and instead focused on average frequencies. Kuramoto's original work came out in 1975, when your average machine had the computational ability of a potato. Even in the two decades that followed, if you wanted to visualize a simulation, the best you could do was a single image. Given such restrictions, it made sense to choose the variables that seemed the most information dense, namely, phase and average frequency. And given that most phase images tend to look like an explosion in a Skittles\textsuperscript{TM} factory, the choice was pretty clear. 

This is not to say that examining average frequencies has not been intensely valuable and interesting. However, we should also keep our minds open to the instantaneous dynamics of our dynamical systems. Given our modern set of tools, it is no longer terribly hard to create real-time, interactive versions of even modestly sized Kuramoto lattices. And as \href{https://www.youtube.com/playlist?list=PL-OZTwOSWQ2A8pa1xCvX07y3_I2aYCos9}{\underline{these additional videos}} show, we can have an awful lot of fun simply watching how things evolve. 

All that being said, we still have a lot of questions and directions to take this problem. Figuring out a good way to analytically estimate $\omega_c$ for the frequency spiral is one of them. We have a numerical estimate of $\omega_c \approx 0.40$, but we do not yet see how to derive this critical value theoretically. 

Likewise, we have numerical evidence that frequency spirals decay (roughly) exponentially in space, but we have not managed to calculate their decay constant directly from the governing equations.

There are a number of related, and potentially more tractable, problems that can perhaps shed light on the subject (for example, continuous analogs of the discrete lattice, or networks with simpler topologies). 

Given the start made here, it also seems plausible that we should be able to apply our results on the canonical frequency spiral to the Kuramoto grid at the (still mysterious) onsets of phase locking and frequency locking discussed in Section 2. Thanks to the work of Leet et al.~\cite{lee10}, we now have a better computational understanding of the dynamics underlying a phase vortex, as well as its effects on its neighborhood. However, we have yet to clearly describe (and ideally, derive) the rules of motion for these spirals, though our presented observations make for a good foundation. It has already been noted that the frequency spiral has a phase vortex at its core, and the motion of these vortices is a major contributor to the establishment of average-frequency clusters. So we may be able to use our understanding of instantaneous dynamics to further our ability to make statements about frequency locking.

Speaking of related problems, the canonical frequency spiral falls into a larger class of problems involving systems of uniform oscillators with a single defector. The situation we spent this paper examining is rather narrow! Pure sine-coupled oscillators are likely to be a non-generic case, given their odd symmetry and their single harmonic character. What happens if these assumptions are relaxed?  In particular, there are very good reasons to believe that even a small degree of non-oddness can qualitatively change how a network of coupled oscillators behaves~\cite{sakaguchi88, ostborn04}. Two-dimensional grids are a specific, strangely difficult topology. There are likely other topologies we could use to inform our understanding of this problem, such as a binary tree. And what about defects that act as pure drivers, forcing the lattice of oscillators without being affected by them in return? Perhaps forced nonlinear lattices could emulate some of the phenomena seen here in a more analytically tractable form. But this is all for future days.

This research was supported by a Sloan Fellowship to Bertrand Ottino-L\"{o}ffler, in the Center for Applied Mathematics in Cornell, as well as by NSF grant DMS-1513179 to Steven Strogatz.

\end{document}